# PROBABILITIES OF RANDOMLY CENTERED SMALL BALLS AND QUANTIZATION IN BANACH SPACES

BY S. DEREICH[1] AND M. A. LIFSHITS[2]

*Technische Universität Berlin and St. Petersburg State University*

We investigate the Gaussian small ball probabilities with random centers, find their deterministic a.s.-equivalents and establish a relation to infinite-dimensional high-resolution quantization.

**1. Introduction.** Consider a centered Gaussian vector $X$ in a separable Banach space $(E, \| \cdot \|)$ with law $\mu$ and reproducing kernel Hilbert space (RKHS) $(H, \| \cdot \|_\mu)$. We let $B$ and $B_\mu$ denote the closed unit balls in $E$ and $H$, respectively. We also use the following notation for shifted balls: $B(x, \varepsilon) := x + \varepsilon B$. The small ball function (SBF) $\varphi$ is defined by

$$\varphi(\varepsilon) := -\log \mu(B(0, \varepsilon)), \qquad \varepsilon > 0.$$

The properties of SBF have been extensively investigated during the last decade; see, for example, [11] and [13]. See also works [6] and [8] on further deep applications of SBF. A complete bibliography on the topic can be found on the website http://www.proba.jussieu.fr/pageperso/smalldev.

Our aim is to study the concentration properties of the r.v.

$$\ell_\varepsilon(\omega) := -\log \mu(B(X(\omega), \varepsilon)), \qquad \varepsilon > 0,$$

the *random small ball function* (RSBF), when $\varepsilon$ is small. We will see that some typical features of the SBF are true as well for the RSBF but the exact asymptotics of the two functions do not coincide.

Beyond structural properties of Gaussian measures, the research is motivated by a close link to so-called random strategies in quantization problems, that we briefly recall now. Let $E$ be a space of objects (images, pictures,

Received February 2004; revised September 2004.

[1]Author's Ph.D. thesis provided main basis for this work.

[2]Supported by Grants RFBR 02-01-00265 and NSh 2258.2003.1.

*AMS 2000 subject classifications.* 60G35, 46F25, 94A15.

*Key words and phrases.* High-resolution quantization, small ball probabilities, small deviations, asymptotic equipartition property.

speech records, etc.) we want to code via a finite codebook. In particular, one can take a finite subset of $E$, say, $(y_i)_{i \leq n}$, as a codebook. In the spirit of Bayesian approach, assume that the subject of coding $X \in E$ is random and its distribution (prior measure $\mu$) is known. Then we can evaluate the quality of a codebook (quantization error) by

$$d(s) = \mathbb{E} \left[ \min_{i=1,\ldots,n} \| X - y_i \|^s \right]^{1/s}.$$

In general it is not feasible to find optimal codebooks under a given constraint on the size $n$ of the codebook. Therefore, recent research focused on the finding of asymptotically good codebooks or on the determination of the (weak or strong) asymptotics of the (theoretically) best achievable coding quality when $n$ tends to infinity, the so-called *high-resolution quantization problem*. It was shown in [5] that these weak asymptotics are in many cases of the same order as the inverse of the small ball function. If the underlying space is a Hilbert space and under a polynomial decay assumption on the eigenvalues of the covariance operator, Luschgy and Pagès [15] proved equivalence of the strong asymptotics to the Shannon distortion rate function. Now using an explicit formula for the distortion rate function based on the eigenvalues, the problem can often be solved explicitly.

In the general high-resolution case, a reasonable codebook can be created by taking independent $\mu$-distributed variables $\{Y_i\}$ (assuming also their independence of $X$). We are thus led to consider the approximation quantities

$$D(r,s) = \mathbb{E} \left[ \min_{i=1,\ldots,\lfloor e^r \rfloor} \| X - Y_i \|^s \right]^{1/s}.$$

The asymptotics of $D(r,s)$, $r \to \infty$, were related to the (standard) small ball function in [5]. Some first properties of the random small ball function and its close relationship to the asymptotics of $D(\cdot, s)$ have been derived in [4]. Whenever the underlying space $E$ is a separable Hilbert space, the RSBF is almost surely equivalent to an invertible deterministic function $\varphi_* : \mathbb{R}_+ \to \mathbb{R}_+$. Moreover, one has

$$D(r,s) \sim \varphi_*^{-1}(r), \qquad r \to \infty,$$

for arbitrary $s > 0$ under certain assumptions on the eigenvalues of the underlying covariance operator. Here and elsewhere we write $f \sim g$ iff $\lim \frac{f}{g} = 1$, while $f \lesssim g$ stands for $\limsup \frac{f}{g} \leq 1$. Finally, $f \approx g$ means

$$0 < \liminf \frac{f}{g} \leq \limsup \frac{f}{g} < \infty.$$

In this article we extend all mentioned results to the Banach space setting. Since the proofs in [4] made strong use of the Hilbertian structure, the new techniques used here differ significantly from those used previously.

The article is arranged as follows. First we prove an almost sure upper bound for $\ell_\varepsilon$ based on the SBF. In Section 3 we find a.s.-equivalence of the RSBF and its median under weak regularity conditions. Some alternative gauge functions for the RSBF are considered in Section 4. In Section 5 a link between the approximation quantity $D$ and the RSBF is established. Finally, in Section 6 the existence of polynomial equivalents for the RSBF is shown in some important particular cases.

## 2. General properties of RSBF.

THEOREM 2.1. *One has*

$$\ell_\varepsilon \lesssim 2\varphi(\varepsilon/2) \qquad as \ \varepsilon \downarrow 0, \ a.s.$$

PROOF. For $n \in \mathbb{N}$, denote $c_n = n$ and $\varepsilon_n = \varphi^{-1}(n^3)$. Let $\Phi$ and $\Upsilon$ denote the distribution function and the tail of the standard normal law. Consider the sets (enlarged balls, in Talagrand's terminology, see [17])

$$A_n = \varepsilon_n B + (c_n + \Upsilon^{-1}(\mu(B(0, \varepsilon_n))))B_\mu.$$

Then, by the isoperimetric inequality (see, e.g., [12], Chapter 11):

$$\mu(A_n) \geq \Phi[c_n + \Upsilon^{-1}(\mu(B(0, \varepsilon_n))) + \Phi^{-1}(\mu(B(0, \varepsilon_n)))] = \Phi(c_n).$$

The tail probabilities of standard normal random variables satisfy

$$(2.1) \qquad \Upsilon(y) \leq \tfrac{1}{2}e^{-y^2/2}, \qquad y \geq 0.$$

Therefore,

$$\sum_{n \in \mathbb{N}} \mu(A_n^c) \leq \sum_{n \in \mathbb{N}} \Upsilon(c_n) < \infty.$$

By the Borel–Cantelli lemma, almost surely all but finitely many events $\{X \in A_n\}$, $n \in \mathbb{N}$, occur.

On the other hand, for every $x \in A_n$ there exists $h \in H$ such that $|h|_\mu \leq c_n + \Upsilon^{-1}(\mu(B(0, \varepsilon_n)))$ and $\|x - h\| \leq \varepsilon_n$; thus, using Borell's shift inequality (see, e.g., [12], page 150), one has

$$(2.2) \qquad \begin{aligned} \mu(B(x, 2\varepsilon_n)) &\geq \mu(B(h, \varepsilon_n)) \geq \exp\left\{-\frac{|h|_\mu^2}{2} - \varphi(\varepsilon_n)\right\} \\ &\geq \exp\left\{-\frac{1}{2}[c_n + \Upsilon^{-1}(\mu(B(0, \varepsilon_n)))]^2 - \varphi(\varepsilon_n)\right\}. \end{aligned}$$

Using the elementary consequence of (2.1)

$$(2.3) \qquad \Upsilon^{-1}(u) \leq \sqrt{-2\log u}, \qquad u \in (0, 1/2],$$

we arrive at

$$-\log \mu(B(x, 2\varepsilon_n)) \leq \frac{1}{2}[c_n + \Upsilon^{-1}(\mu(B(0, \varepsilon_n)))]^2 + \varphi(\varepsilon_n)$$

$$\leq \frac{1}{2}[c_n + \sqrt{2\varphi(\varepsilon_n)}]^2 + \varphi(\varepsilon_n)$$

$$= \frac{c_n^2}{2} + c_n\sqrt{2\varphi(\varepsilon_n)} + 2\varphi(\varepsilon_n).$$

Note that $c_n^2 = o(\varphi(\varepsilon_n))$ as $n \to \infty$ and, therefore,

$$\sup_{x \in A_n} -\log \mu(B(x, 2\varepsilon_n)) \lesssim 2\varphi(\varepsilon_n), \qquad n \to \infty.$$

Since $\lim_{n\to\infty} \varphi(\varepsilon_{n+1})/\varphi(\varepsilon_n) = 1$ and the small ball probabilities are monotone, our theorem is proved. $\square$

REMARK 2.2. The previous theorem and Anderson's inequality (see, e.g., [12], page 135) imply that the random small ball function $\ell_\varepsilon$ is asymptotically enclosed between two deterministic functions, that is,

(2.4)                    $\varphi(\varepsilon) \leq \ell_\varepsilon \lesssim 2\varphi(\varepsilon/2), \qquad \varepsilon \downarrow 0,$  a.s.

Suppose now that there exists $\nu < \infty$ such that

(2.5)                    $\varphi(\varepsilon) \leq \nu\varphi(2\varepsilon)$

for sufficiently small $\varepsilon > 0$. Then the RSBF function is of the same order as the small ball function and we have

$$\varphi(\varepsilon) \leq \ell_\varepsilon \lesssim 2\nu\varphi(\varepsilon), \qquad \varepsilon \downarrow 0,$$  a.s.

A better asymptotic lower bound will be presented in Corollary 4.4 below.

REMARK 2.3. One can find alternative estimates for probabilities of enlarged balls $\varepsilon B + r B_\mu$ in [17]. These estimates proved to be more efficient than the isoperimetric inequality in the work concerning Strassen's functional law of the iterated logarithm, where they yield the correct convergence rate. Surprisingly, in the range of parameters $\varepsilon, r$ considered in our work, the estimates from [17] provide worse results than the isoperimetric inequality.

**3. Equivalence to a deterministic function.** The main objective of this section is to prove concentration inequalities for the random variables $\ell_\varepsilon$ as $\varepsilon \downarrow 0$. In the main theorem, we will find equivalence of random small ball probabilities to a deterministic function under weak assumptions.

It is well known that concentration phenomena occur for $H$-Lipschitz functionals. We will show, by using a result of Kuelbs and Li [6], that the function $\log \mu(B(\cdot, \varepsilon))$ is $H$-Lipschitz on a set of probability "almost 1," and the corresponding Lipschitz constant will be controlled.

3.1. *Large set of good points.* Let us fix $\varepsilon > 0$ and choose $M = M(\varepsilon) = 3\sqrt{\varphi(\varepsilon)}$. Introduce again an enlarged ball

$$V_\varepsilon := \varepsilon B + M B_\mu.$$

Let us start by showing that $V_\varepsilon$ is large enough. Indeed, by the isoperimetric inequality and (2.3):

$$
\begin{aligned}
\mu(V_\varepsilon^c) &\leq \Upsilon(\Phi^{-1}(\mu(\varepsilon B)) + M) \\
&= \Upsilon(-\Upsilon^{-1}(\mu(\varepsilon B)) + M) \\
&\leq \Upsilon(-\sqrt{-2\log\mu(\varepsilon B)} + M) \\
&= \Upsilon(-\sqrt{2\varphi(\varepsilon)} + M) \\
&\leq \exp(-\varphi(\varepsilon)).
\end{aligned}
\tag{3.1}
$$

We also observe that the small ball probabilities are uniformly bounded from below on $V_\varepsilon$. Indeed, for each $x \in V_\varepsilon$, there exists $h \in M B_\mu \cap B(x, \varepsilon)$. Hence, $B(x, 2\varepsilon) \supset B(h, \varepsilon)$, and we obtain, similarly to (2.2),

$$
\begin{aligned}
\log\mu(B(x, 2\varepsilon)) &\geq \log\mu(B(h, \varepsilon)) \\
&\geq \log(\exp(-|h|_\mu^2/2)\mu(B(0, \varepsilon))) \\
&\geq -M^2/2 - \varphi(\varepsilon) = -5.5\varphi(\varepsilon).
\end{aligned}
\tag{3.2}
$$

3.2. *Estimate of the Lipschitz constant.* In this section, we consider the $H$-Lipschitz property of the function $\Psi(\cdot) := \log\mu(B(\cdot, 2\varepsilon))$ on $V_\varepsilon$.

PROPOSITION 3.1. *Let $\varepsilon > 0$ be so small that*

$$\varphi(2\varepsilon) \geq -\log\Phi(-3).$$

*Let $h \in H$ and $x, x + h \in V_\varepsilon$. Then*

$$|\Psi(x + h) - \Psi(x)| \leq 8\sqrt{\varphi(\varepsilon)}|h|_\mu. \tag{3.3}$$

PROOF. Since $V_\varepsilon$ is convex, without loss of generality we may and do assume that $|h|_\mu \leq 1$. Since $B(x + h, 2\varepsilon) = B(x, 2\varepsilon) + h$, we can use the estimate from [7] which states that for an arbitrary Gaussian measure $\mu$, a measurable set $A$ and an element $h$ of the RKHS $H$, one has

$$\Phi(\Phi^{-1}(\mu(A)) - |h|_\mu) \leq \mu(A + h) \leq \Phi(\Phi^{-1}(\mu(A)) + |h|_\mu). \tag{3.4}$$

Thus,

$$\mu(B(x + h, 2\varepsilon)) \leq \Phi(\theta + |h|_\mu),$$

where $\theta = \Phi^{-1}(\mu(B(x, 2\varepsilon)))$, and we obtain

$$
\begin{aligned}
(3.5) \qquad \Delta &:= \Psi(x + h) - \Psi(x) \\
&= \log \mu(B(x + h, 2\varepsilon)) - \log \mu(B(x, 2\varepsilon)) \\
&\leq \log \Phi(\theta + |h|_\mu) - \log \Phi(\theta).
\end{aligned}
$$

Under our assumptions it is true that

$$
\mu(B(0, 2\varepsilon)) = \exp(-\varphi(2\varepsilon)) \leq \Phi(-3).
$$

Therefore,

$$
\Phi(\theta) = \mu(B(x, 2\varepsilon)) \leq \mu(B(0, 2\varepsilon)) \leq \Phi(-3),
$$

which shows that $\theta \leq -3$. It follows from $|h|_\mu \leq 1$ that $\theta \leq \theta + |h|_\mu \leq -2$. Using the elementary inequality

$$
0 \leq (\log \Phi)'(r) \leq 2|r|, \qquad r \leq -2,
$$

we obtain

$$
\log \Phi(\theta + |h|_\mu) - \log \Phi(\theta) \leq 2 \int_\theta^{\theta + |h|_\mu} |r| \, dr = 2|\theta||h|_\mu - |h|_\mu^2 \leq 2|\theta||h|_\mu.
$$

Now note that due to (3.2)

$$
\Phi(\theta) = \mu(B(x, 2\varepsilon)) \geq \exp(-5.5\varphi(\varepsilon))
$$

and, hence by (2.3),

$$
|\theta| = \Upsilon^{-1}(\mu(B(x, \varepsilon))) \leq \sqrt{11\varphi(\varepsilon)}.
$$

Altogether, we obtain

$$
\Delta \leq 8\sqrt{\varphi(\varepsilon)}|h|_\mu.
$$

To derive the converse bound, we use that the situation is symmetric. Namely, take $\tilde{x} = x + h$ and $\tilde{h} = -h$. Then we have $\tilde{x}, \tilde{x} + \tilde{h} \in V_\varepsilon$ and the arguments from above imply that

$$
\begin{aligned}
-\Delta = \Psi(x) - \Psi(x + h) &= \Psi(\tilde{x} + \tilde{h}) - \Psi(\tilde{x}) \\
&\leq 8\sqrt{\varphi(\varepsilon)}|h|_\mu. \qquad\qquad \square
\end{aligned}
$$

3.3. *Concentration and convergence.* We are now in a position to prove our main result on the deterministic equivalent for the RSBF.

THEOREM 3.2. *Assume that for all $\varepsilon > 0$ small enough it is true that*

$$(3.6) \qquad \varphi(\varepsilon) \leq \nu\varphi(2\varepsilon)$$

*for some $\nu < \infty$. Let $m_\varepsilon$ be a median of $\ell_\varepsilon$. Then*

$$\lim_{\varepsilon \downarrow 0} \frac{\ell_\varepsilon}{m_\varepsilon} = 1 \qquad \text{almost surely.}$$

PROOF. Define $r_\varepsilon$ from the equation

$$\Phi(r_\varepsilon) = \tfrac{1}{2} + \mu(V_\varepsilon^c).$$

It follows from (3.1) that $\lim_{\varepsilon \downarrow 0} r_\varepsilon = 0$. By (3.3) and the concentration principle for $H$-Lipschitz functionals (see, e.g., [9], page 210) we have, for any $r > r_\varepsilon$ that,

$$\mathbb{P}(|\ell_{2\varepsilon} - m_{2\varepsilon}| \geq 8\sqrt{\varphi(\varepsilon)}r) \leq \mu(V_\varepsilon^c) + \exp(-(r - r_\varepsilon)^2/2).$$

Let us fix $\delta > 0$ and let

$$r = \frac{\delta\varphi(2\varepsilon)}{8\sqrt{\varphi(\varepsilon)}}.$$

Then using (3.1), we obtain that

$$\mathbb{P}(|\ell_{2\varepsilon} - m_{2\varepsilon}| \geq \delta\varphi(2\varepsilon)) \leq \exp(-\varphi(\varepsilon)) + \exp\left(-\frac{\delta^2\varphi(2\varepsilon)^2(1 + o(1))}{2 \cdot 8^2\varphi(\varepsilon)}\right).$$

Due to (3.6), it holds for $\varepsilon > 0$ sufficiently small

$$\begin{aligned} \mathbb{P}(|\ell_{2\varepsilon} - m_{2\varepsilon}| \geq \delta\varphi(2\varepsilon)) &\leq \exp(-\varphi(\varepsilon)) + \exp\left(-\frac{\delta^2\varphi(2\varepsilon)}{3 \cdot 8^2\nu}\right) \\ &\leq 2\exp\left(-\min\left\{1; \frac{\delta^2}{3 \cdot 8^2\nu}\right\}\varphi(2\varepsilon)\right) \\ &=: 2\exp\{-\nu_1\varphi(2\varepsilon)\}. \end{aligned}$$

By switching from $2\varepsilon$ to $\varepsilon$, we get

$$(3.7) \qquad \mathbb{P}(|\ell_\varepsilon - m_\varepsilon| \geq \delta\varphi(\varepsilon)) \leq 2\exp\{-\nu_1\varphi(\varepsilon)\}.$$

Next recall that due to Remark 2.2, $\ell_\varepsilon$ is asymptotically a.s. enclosed by the two functions $\varphi(\varepsilon)$ and $2\nu\varphi(\varepsilon)$. Hence, it holds for $\varepsilon > 0$ sufficiently small

$$(3.8) \qquad \varphi(\varepsilon) \leq m_\varepsilon \leq 3\nu\varphi(\varepsilon).$$

Now consider for $n \in \mathbb{N}$ the set $T_n := \{\varepsilon > 0 : n \leq m_\varepsilon < n + 1\}$. For sufficiently large $n \in \mathbb{N}$ and $\varepsilon \in T_n$ we have, by using (3.7) and (3.8),

$$\begin{aligned}
\mathbb{P}(\ell_\varepsilon > (1+\delta)(n+1)) &\leq \mathbb{P}(\ell_\varepsilon > (1+\delta)m_\varepsilon) \\
&\leq \mathbb{P}(\ell_\varepsilon - m_\varepsilon > \delta m_\varepsilon) \\
&\leq 2 \exp\{-\nu_1 \varphi(\varepsilon)\} \\
&\leq 2 \exp\left\{-\frac{\nu_1}{3\nu} m_\varepsilon\right\} \\
&\leq 2 \exp\left\{-\frac{\nu_1}{3\nu} n\right\}
\end{aligned}$$

so that

$$\mathbb{P}\left(\sup_{T_n} \ell_\varepsilon > (1+\delta)(n+1)\right) \leq 2 \exp\left\{-\frac{\nu_1}{3\nu} n\right\}.$$

By the Borel–Cantelli lemma we eventually have, for all large $n$ and all $\varepsilon \in T_n$,

$$\ell_\varepsilon \leq (1+\delta)(n+1) \leq (1+\delta)\frac{n+1}{n} m_\varepsilon,$$

and, since $\delta > 0$ can be chosen arbitrarily small, it follows that

$$\limsup_{\varepsilon \downarrow 0} \frac{\ell_\varepsilon}{m_\varepsilon} \leq 1 \qquad \text{a.s.}$$

The inverse bound can be obtained in the same way.  □

In the case where Theorem 3.2 is not applicable, we still can show:

PROPOSITION 3.3.   *For any continuous function $\psi : \mathbb{R}_+ \to \mathbb{R}_+$ such that $\lim_{\varepsilon \downarrow 0} \psi(\varepsilon) = \infty$, there exist constants $c_\psi, C_\psi \in [0, \infty]$ such that*

$$\liminf_{\varepsilon \downarrow 0} \frac{\ell_\varepsilon}{\psi(\varepsilon)} = c_\psi$$

*and*

$$\limsup_{\varepsilon \downarrow 0} \frac{\ell_\varepsilon}{\psi(\varepsilon)} = C_\psi$$

*almost surely.*

PROOF.   Let $E^*$ denote the topological dual of $E$ and $C_\mu : E^* \to E$ the covariance operator of $\mu$. Let $x \in E$ and $h = C_\mu(z)$ for some $z \in E^*$. The Cameron–Martin formula (see, e.g., [12], page 107) gives

$$\mu(B(x-h, \varepsilon)) = \int_{B(x,\varepsilon)} \exp\{z(y) - \tfrac{1}{2}\|z\|^2_{L_2(\mu)}\} d\mu(y).$$

Since $z : E \to \mathbb{R}$ is a continuous function, it holds

$$\mu(B(x - h, \varepsilon)) \sim \exp\{z(x) - \tfrac{1}{2}\|z\|_{L_2(\mu)}^2\} \mu(B(x, \varepsilon))$$

as $\varepsilon \downarrow 0$. In particular,

$$-\log \mu(B(x - h, \varepsilon)) \sim -\log \mu(B(x, \varepsilon)), \qquad \varepsilon \downarrow 0.$$

Therefore, for any $s \geq 0$, the set

$$A_s = \left\{ x \in E : \liminf_{\varepsilon \downarrow 0} \frac{-\log \mu(B(x, \varepsilon))}{\psi(\varepsilon)} \leq s \right\}$$

is invariant under an arbitrary shift $h \in C_\mu(E^*)$. Since $\psi$ is continuous, the set $A_s$ is measurable. Moreover, by the zero–one law for Gaussian measures (see [1], Theorem 2.5.2), the set $A_s$ has $\mu$-measure 0 or 1. The first assertion follows. The second one may be proved analogously. $\square$

## 4. Gauge functions.

In this section we suppose that the regularity condition (3.6) applies. By Theorem 3.2 one has

$$(4.1) \qquad \qquad \ell_\varepsilon \sim \varphi_*(\varepsilon) \qquad \text{as } \varepsilon \downarrow 0, \text{ a.s.}$$

for $\varphi_*(\varepsilon)$, $\varepsilon > 0$, equal to the median $m_\varepsilon$ of $\ell_\varepsilon$. In the sequel, we study alternative representations for $\varphi_*$. We will need the following lemma.

LEMMA 4.1. *Let $Z$ denote a standard normal r.v. For any $p \geq 1$ and $\varepsilon > 0$ with $\mu(\varepsilon B) \leq 1/2$, one has*

$$\|\ell_{2\varepsilon}\|_{L^p(\mathbb{P})} \leq \varphi(\varepsilon) + \tfrac{1}{2}(\sqrt{2\varphi(\varepsilon)} + \|Z\|_{L^{2p}(\mathbb{P})})^2.$$

PROOF. The proof is similar to that of Theorem 2.1. We fix $\varepsilon > 0$ with $\mu(\varepsilon B) \leq 1/2$ and let

$$A_t = \varepsilon B + (t + \Upsilon^{-1}(\mu(\varepsilon B))) B_\mu, \qquad t \geq 0.$$

By the isoperimetric inequality one has

$$(4.2) \qquad \mu(A_t) \geq \Phi[t + \Upsilon^{-1}(\mu(\varepsilon B)) + \Phi^{-1}(\mu(\varepsilon B))] = \Phi(t)$$

for any $t \geq 0$. As in the proof of Theorem 2.1, we obtain for $x \in A_t$,

$$\mu(B(x, 2\varepsilon)) \geq \exp\{-\tfrac{1}{2}[t + \Upsilon^{-1}(\mu(\varepsilon B))]^2 - \varphi(\varepsilon)\}$$

and inequality (2.3) yields

$$-\log \mu(B(x, 2\varepsilon)) \leq \tfrac{1}{2}[t + \Upsilon^{-1}(\mu(\varepsilon B))]^2 + \varphi(\varepsilon)$$

$$\leq \tfrac{1}{2}[t + \sqrt{2\varphi(\varepsilon)}]^2 + \varphi(\varepsilon).$$

Combining this estimate with (4.2) gives

$$\mathbb{P}(\ell_{2\varepsilon} > \tfrac{1}{2}[t + \sqrt{2\varphi(\varepsilon)}]^2 + \varphi(\varepsilon)) \le \Upsilon(t)$$

for all $t \ge 0$. Hence, with $Z^+ = Z \vee 0$ it follows that

$$\|\ell_{2\varepsilon}\|_{L^p(\mathbb{P})} \le \mathbb{E}[(\tfrac{1}{2}[Z^+ + \sqrt{2\varphi(\varepsilon)}]^2 + \varphi(\varepsilon))^p]^{1/p}.$$

Applying the triangle inequality twice yields

$$\|\ell_{2\varepsilon}\|_{L^p(\mathbb{P})} \le \tfrac{1}{2}\mathbb{E}[(Z^+ + \sqrt{2\varphi(\varepsilon)})^{2p}]^{1/p} + \varphi(\varepsilon)$$

$$\le \tfrac{1}{2}(\mathbb{E}[(Z^+)^{2p}]^{1/2p} + \sqrt{2\varphi(\varepsilon)})^2 + \varphi(\varepsilon)$$

and the assertion follows. $\square$

THEOREM 4.2.  *For $\varphi_*$ satisfying* (4.1), *we have*

$$\lim_{\varepsilon \downarrow 0} \frac{\ell_\varepsilon}{\varphi_*(\varepsilon)} = 1$$

*in $L^p(\mathbb{P})$ for any $p \ge 1$. In particular,*

$$\varphi_*(\varepsilon) \sim \|\ell_\varepsilon\|_{L^p(\mathbb{P})} \qquad as\ \varepsilon \downarrow 0.$$

PROOF.  Fix $\eta \in (0,1)$ and let

$$\mathcal{T}(\varepsilon) = \left\{ x \in E : \left| \frac{-\log \mu(B(x,\varepsilon))}{\varphi_*(\varepsilon)} - 1 \right| \le \eta \right\}.$$

Recall that $X = X(\omega)$ denotes the $\mu$-distributed center of the random ball. One has

$$\left\| \frac{\ell_\varepsilon}{\varphi_*(\varepsilon)} - 1 \right\|_{L^p(\mathbb{P})} \le \left\| \mathbb{1}_{\mathcal{T}(\varepsilon)}(X)\left( \frac{\ell_\varepsilon}{\varphi_*(\varepsilon)} - 1 \right) \right\|_{L^p(\mathbb{P})}$$

$$+ \left\| \mathbb{1}_{\mathcal{T}(\varepsilon)^c}(X)\left( \frac{\ell_\varepsilon}{\varphi_*(\varepsilon)} - 1 \right) \right\|_{L^p(\mathbb{P})}$$

$$\le \left\| \mathbb{1}_{\mathcal{T}(\varepsilon)}(X)\left( \frac{\ell_\varepsilon}{\varphi_*(\varepsilon)} - 1 \right) \right\|_{L^p(\mathbb{P})}$$

$$+ \left\| \mathbb{1}_{\mathcal{T}(\varepsilon)^c}(X)\frac{\ell_\varepsilon}{\varphi_*(\varepsilon)} \right\|_{L^p(\mathbb{P})} + \left\| \mathbb{1}_{\mathcal{T}(\varepsilon)^c}(X) \right\|_{L^p(\mathbb{P})}$$

$$=: I_1(\varepsilon) + I_2(\varepsilon) + I_3(\varepsilon).$$

Clearly, $I_1(\varepsilon) \le \eta$. Using the Cauchy–Schwarz inequality, we estimate the second term by

$$I_2(\varepsilon) = \frac{1}{\varphi_*(\varepsilon)} \|\mathbb{1}_{\mathcal{T}(\varepsilon)^c}(X)\ell_\varepsilon\|_{L^p(\mathbb{P})}$$

$$\le \frac{1}{\varphi_*(\varepsilon)} \mu(\mathcal{T}(\varepsilon)^c)^{1/2p} \|\ell_\varepsilon\|_{L^{2p}(\mathbb{P})}.$$

By (3.6) and the previous lemma, $\|\ell_\varepsilon\|_{L^{2p}(\mathbb{P})} \lesssim 2\varphi(\varepsilon/2) \lesssim 2\nu\varphi(\varepsilon)$ as $\varepsilon \downarrow 0$. Due to Anderson's inequality, one has $\varphi_*(\varepsilon) \gtrsim \varphi(\varepsilon)$ ($\varepsilon \downarrow 0$). On the other hand, $\lim_{\varepsilon\downarrow 0} \mu(\mathcal{T}(\varepsilon)^c) = 0$ by assumption (4.1). Hence,

$$\lim_{\varepsilon\downarrow 0} I_2(\varepsilon) = 0.$$

Furthermore, $\lim_{\varepsilon\downarrow 0} I_3(\varepsilon) = 0$. Putting all three estimates together gives

$$\left\| \frac{\ell_\varepsilon}{\varphi_*(\varepsilon)} - 1 \right\|_{L^p(\mathbb{P})} \lesssim \eta, \qquad \varepsilon \downarrow 0.$$

Since $\eta \in (0,1)$ was arbitrary, the proof is complete. □

REMARK 4.3. As a consequence of the above theorem one can replace the median $m_\varepsilon$ by $\mathbb{E}[\ell_\varepsilon]$ in Theorem 3.2. By the well-known fact that small ball functions are convex, it is easy to deduce that the function $\mathbb{R}_+ \to \mathbb{R}_+$, $\varepsilon \mapsto \mathbb{E}[\ell_\varepsilon]$ is convex, one-to-one and onto.

COROLLARY 4.4. *It is true that*

$$\varphi(\varepsilon/\sqrt{2}) \lesssim \varphi_*(\varepsilon) \lesssim 2\varphi(\varepsilon/2), \qquad \varepsilon \downarrow 0.$$

PROOF. The asymptotic upper bound follows from Theorem 2.1. It remains to prove the lower bound. Due to the previous remark we can prove the statement for $\varphi_*(\varepsilon) = \mathbb{E}[\ell_\varepsilon]$. Denote by $\tilde{X}$ a $\mu$-distributed r.v. that is independent of $X$. One has for $\varepsilon > 0$,

$$-\mathbb{E}[\ell_\varepsilon] = \mathbb{E}[\log \mu(B(X, \varepsilon))]$$

$$= \mathbb{E}[\log \mathbb{P}(\|X - \tilde{X}\| \le \varepsilon | X)]$$

$$\le \log \mathbb{E}[\mathbb{P}(\|X - \tilde{X}\| \le \varepsilon | X)]$$

$$= \log \mathbb{P}(\|X - \tilde{X}\| \le \varepsilon),$$

where the inequality follows from Jensen's inequality. Note that $X - \tilde{X}$ and $\sqrt{2}X$ are both centered Gaussian vectors with the same covariance operator. Therefore, $\mathcal{L}(\sqrt{2}X) = \mathcal{L}(X - \tilde{X})$, which shows that

$$\mathbb{E}[\ell_\varepsilon] \ge -\log \mu(B(0, \varepsilon/\sqrt{2})) = \varphi(\varepsilon/\sqrt{2}). \qquad\qquad □$$

REMARK 4.5. If there exists $\tilde{\nu} > 1$ such that

$$\varphi(\varepsilon) \geq \tilde{\nu}\varphi(\sqrt{2}\varepsilon)$$

for all sufficiently small $\varepsilon > 0$, then the strong asymptotics of $\varphi$ and $\varphi_*$ differ.

**5. Equivalence of random small ball probabilities and random quantization.** In this section, we relate the asymptotics of $D(\cdot, s)$ to the RSBF. Recall that

$$D(r, s) = \mathbb{E}\left[\min_{i=1,\ldots,\lfloor e^r \rfloor} \|X - Y_i\|^s\right]^{1/s},$$

where $\{Y_i\}_{i \in \mathbb{N}}$ is a sequence of independent (of $X$ as well) $\mu$-distributed r.v.'s in $E$. In terms of information theory, Theorem 3.2 can be interpreted as the *asymptotic equipartition property* (AEP) corresponding to the random quantization problem. For more details on AEPs and their connections to coding theory we refer the reader to [2].

In the sequel, we assume the existence of a convex function $\varphi_* : \mathbb{R}_+ \to \mathbb{R}_+$ that is one-to-one, onto and satisfies

(5.1)        $\varphi_*(\varepsilon) \sim -\log \mu(B(X, \varepsilon))$        as $\varepsilon \downarrow 0$, in probability.

By the preceding considerations $\varphi_*(\varepsilon) := \mathbb{E}[\ell_\varepsilon]$ is an appropriate choice if (3.6) applies.

THEOREM 5.1. *If there exists $\tilde{\nu} > 1$ such that for all $\varepsilon > 0$ sufficiently small*

(5.2)                          $\varphi(\varepsilon) \geq \tilde{\nu}\varphi(2\varepsilon),$

*then*

$$D(r, s) \sim \varphi_*^{-1}(r), \qquad r \to \infty.$$

We will need a couple of elementary results:

LEMMA 5.2. *Let $f : \mathbb{R}_+ \to \mathbb{R}_+$ be a decreasing convex function satisfying*

$$f(2r) \geq \nu f(r)$$

*for all $r$ sufficiently large. Then, for any function $\Delta : \mathbb{R}_+ \to \mathbb{R}$, with $\Delta(r) = o(r)$ $(r \to \infty)$, one has*

(5.3)                    $f(r + \Delta r) \sim f(r)$        as $r \to \infty.$

PROOF. Convexity yields that for all $\delta, r > 0$ it is true that

$$f((1 + \delta)r) \geq \left(1 - \frac{2\delta}{\nu}\right) f(r).$$

The further necessary estimates are trivial.  □

LEMMA 5.3. *Assumption* (5.2) *implies that there exists a constant* $\nu_1$ *such that*

$$\varphi_*^{-1}(2r) \geq \nu_1 \varphi_*^{-1}(r) \quad and \quad \varphi^{-1}(2r) \geq \nu_1 \varphi^{-1}(r) \tag{5.4}$$

*for sufficiently large* $r \geq 0$.

PROOF. Choose $\kappa \in \mathbb{N}$ such that $6/\tilde{\nu}^\kappa < 1$. By assumption (5.1) and Theorem 2.1, one has for $\varepsilon > 0$ sufficiently small

$$\varphi_*(\varepsilon) \leq 3\varphi(\varepsilon/2) \leq \frac{3}{\tilde{\nu}^\kappa}\varphi(\varepsilon/2^{\kappa+1}) \leq \frac{1}{2}\varphi_*(\varepsilon/2^{\kappa+1}). \tag{5.5}$$

Consequently,

$$\varphi_*^{-1}(2\varphi_*(\varepsilon)) \geq \frac{\varepsilon}{2^{\kappa+1}}.$$

Choosing $\varepsilon = \varphi_*^{-1}(r)$ and assuming that $r$ is sufficiently large, we obtain

$$\varphi_*^{-1}(2r) \geq \frac{1}{2^{\kappa+1}}\varphi_*^{-1}(r).$$

In a similar way, the equation

$$\varphi(\varepsilon) \leq \frac{1}{\tilde{\nu}^\kappa}\varphi(\varepsilon/2^\kappa) \leq \frac{1}{2}\varphi(\varepsilon/2^{\kappa+1})$$

can be used to derive the second assertion of the lemma. □

PROPOSITION 5.4. *Let* $\kappa \in (0, 1)$ *and* $\delta := \frac{1}{4}\min\{1, \nu_1\kappa\}$, *where* $\nu_1$ *satisfies* (5.4). *For* $r \geq 0$ *define the sets*

$$\mathcal{X}_1(r) = \{x \in E : -\log\mu(B(x, (1+\kappa)\varphi_*^{-1}(r))) \leq (1-\delta)r\}$$

*and*

$$\mathcal{X}_2(r) = \{x \in E : -\log\mu(B(x, (1-\kappa)\varphi_*^{-1}(r))) \geq (1+\delta)r\}.$$

*Then*

$$\lim_{r \to \infty}\mu(\mathcal{X}_1(r)) = \lim_{r \to \infty}\mu(\mathcal{X}_2(r)) = 1.$$

PROOF. Making use of the convexity of $\varphi_*^{-1}$ and equation (5.4), one has

$$\varphi_*^{-1}(r - 2\delta r) - \varphi_*^{-1}(r) \leq \frac{2\delta r}{r/2}(\varphi_*^{-1}(r/2) - \varphi_*^{-1}(r)) \leq \frac{4\delta}{\nu_1}\varphi_*^{-1}(r)$$

for $r$ large. Therefore, there exists $r_0 \geq 0$ such that

$$\varphi_*^{-1}(r - 2\delta r) \leq (1+\kappa)\varphi_*^{-1}(r)$$

for all $r \geq r_0$. Consequently, the set $\mathcal{X}_1(r)$ satisfies for $r \geq r_0$

$$\mathcal{X}_1(r) \supset \left\{ x \in E : -\log \mu(B(x, \varphi_*^{-1}((1-2\delta)r))) \leq \frac{1-\delta}{1-2\delta}(1-2\delta)r \right\}.$$

Since $\varphi_*^{-1}((1-2\delta)r)$ converges to 0 and $(1-\delta)/(1-2\delta) > 1$, it holds by assumption (5.1) that

$$\lim_{r \to \infty} \mu(\mathcal{X}_1(r)) = 1.$$

The proof of the second assertion is similar: one has for $r \geq r_0$

$$\varphi_*^{-1}(r) - \varphi_*^{-1}(r+2\delta r) \leq \varphi_*^{-1}(r-2\delta r) - \varphi_*^{-1}(r) \leq \kappa \varphi_*^{-1}(r),$$

where the first inequality is a consequence of the convexity of $\varphi_*^{-1}$. Hence, $\varphi_*^{-1}(r+2\delta r) \geq (1-\kappa)\varphi_*^{-1}(r)$ for $r \geq r_0$ and it follows

$$\mathcal{X}_2(r) \supset \left\{ x \in E : -\log \mu(B(x, \varphi_*^{-1}(r+2\delta r))) \geq \frac{1+\delta}{1+2\delta}(1+2\delta)r \right\}.$$

Finally, assumption (5.1) yields

$$\lim_{r \to \infty} \mu(\mathcal{X}_2(r)) = 1. \qquad \square$$

PROPOSITION 5.5.  *Let $\kappa \in (0,1)$. For $r \geq 0$ consider*

$$Z(r) := \min_{i=1,\dots,\lfloor e^r \rfloor} \|X - Y_i\|$$

*and the event*

$$\mathcal{T}(r) := \{Z(r) \in [(1-\kappa)\varphi_*^{-1}(r), (1+\kappa)\varphi_*^{-1}(r)]\}.$$

*Under the assumptions of Theorem 5.1 one has*

$$\lim_{r \to \infty} \mathbb{P}(\mathcal{T}(r)) = 1.$$

PROOF.  Recall that $\varphi_*^{-1}$ is convex. In view of Lemma 5.2, it suffices to consider $r \in I := \{\log j : j \in \mathbb{N}\}$, that is, the values $r$ for which $e^r$ is an integer. By Proposition 5.4, one has

$$\lim_{r \to \infty} \mu(\mathcal{X}_1(r) \cap \mathcal{X}_2(r)) = 1$$

with $\mathcal{X}_1(r)$ and $\mathcal{X}_2(r)$ as in the proposition. Moreover, for $r \in I$ and $X \in \mathcal{X}_1(r)$, one has

$$\mathbb{P}(Z(r) > (1+\kappa)\varphi_*^{-1}(r)|X) = (1 - \mu(B(X, (1+\kappa)r)))^{e^r}$$

$$\leq (1 - e^{-r+\delta r})^{e^r} = \left(1 - \frac{e^{\delta r}}{e^r}\right)^{e^r}$$

$$\leq \exp\{-e^{\delta r}\} \to 0, \qquad r \to \infty.$$

On the other hand, for $X \in \mathcal{X}_2(r)$, $r \in I$, it holds

$$\mathbb{P}(Z(r) > (1-\kappa)\varphi_*^{-1}(r)|X) = (1 - \mu(B(X,(1-\kappa)r)))^{e^r}$$

$$\geq (1 - e^{-r-\delta r})^{e^r} = \left(1 - \frac{e^{-\delta r}}{e^r}\right)^{e^r} \to 1$$

as $r \to \infty$. Hence, the events $\mathcal{T}(r)$, $r \geq 0$, satisfy $\lim_{r\to\infty} \mathbb{P}(\mathcal{T}(r)) = 1$. $\quad\square$

For the proof of Theorem 5.1, we use a consequence of Theorem 2.1 of [5] (see also [3], Theorem 3.1.2).

THEOREM 5.6. *Suppose there exists $\nu_1 > 0$ such that*

$$\varphi^{-1}(2r) \geq \nu_1 \varphi^{-1}(r)$$

*for all sufficiently large $r \geq 0$. Then one has for arbitrary $s > 0$*

$$D(r,s) \lesssim 2\varphi^{-1}(r/2)$$

*as $r \to \infty$.*

PROOF OF THEOREM 5.1. Fix $s > 0$. First we prove

$$D(r,s) \lesssim \varphi_*^{-1}(r), \qquad r \to \infty.$$

Fix $\kappa \in (0,1)$ and let $\mathcal{T}(r)$ and $Z(r)$ be as in the previous proposition. Now

$$\mathbb{E}[Z(r)^s] \leq \mathbb{E}[\mathbb{1}_{\mathcal{T}(r)}(1+\kappa)^s\varphi_*^{-1}(r)^s] + \mathbb{E}[\mathbb{1}_{\mathcal{T}(r)^c}Z(r)^s] =: I_1(r) + I_2(r).$$

One has $I_1(r) \leq (1+\kappa)^s \varphi_*^{-1}(r)^s$. Moreover, the Cauchy–Schwarz inequality gives

$$I_2(r) \leq \mathbb{P}(\mathcal{T}(r)^c)^{1/2} \mathbb{E}[Z(r)^{2s}]^{1/2}.$$

As a consequence of Lemma 5.3 and assumption (5.2), there exists a constant $\nu_1 > 0$ such that

$$\varphi^{-1}(2r) \geq \nu_1 \varphi^{-1}(r)$$

for large $r \geq 0$. Thus, Theorem 5.6 is applicable and

$$\mathbb{E}[Z(r)^{2s}]^{1/2} \lesssim 2^s \varphi^{-1}(r/2)^s \leq \frac{2^s}{\nu_1^s} \varphi^{-1}(r)^s$$

as $r \to \infty$. By the previous proposition, $\lim_{r\to\infty} \mathbb{P}(\mathcal{T}(r)^c) = 0$. Consequently,

$$I_2(r) = o(\varphi_*^{-1}(r)^s), \qquad r \to \infty$$

and

$$\mathbb{E}[Z(r)^s]^{1/s} \lesssim (1+\kappa)\varphi_*^{-1}(r), \qquad r \to \infty.$$

Since $\kappa \in (0, 1)$ was chosen arbitrarily, it follows that

$$D(r, s) = \mathbb{E}[Z(r)^s]^{1/s} \lesssim \varphi_*^{-1}(r), \qquad r \to \infty.$$

The converse inequality is obvious, since for fixed $\kappa \in (0, 1)$ and $\mathcal{T}(r)$ as above one has

$$\mathbb{E}[Z(r)^s]^{1/s} \geq \mathbb{E}[\mathbb{1}_{\mathcal{T}(r)} Z(r)^s]^{1/s} \geq \mathbb{P}(\mathcal{T}(r))^{1/s}(1 - \kappa)\varphi_*^{-1}(r)$$
$$\gtrsim (1 - \kappa)\varphi_*^{-1}(r), \qquad r \to \infty. \qquad \square$$

**6. Polynomial equivalents for the RSBF.** In the sequel $\mu$ denotes Wiener measure on the canonical Wiener space $\mathbb{C}[0, \infty) := \mathbb{C}([0, \infty), \mathbb{R}^d)$ for some fixed $d \in \mathbb{N}$. Moreover, $X = W$ denotes a $\mathbb{C}[0, \infty)$-valued random variable that is a Wiener process under the standard measure $\mathbb{P}$. We will sometimes use the alternative measures $\mathbb{P}^x$ and $\mathbb{P}^{x,t}$ ($x \in \mathbb{R}^d$, $t \geq 0$). Under these measures $W$ is a Wiener process starting in $x$ at time 0 or at time $t$, respectively. The corresponding expectations are denoted by $\mathbb{E}^x$ and $\mathbb{E}^{x,t}$.

Recall that for many underlying Banach spaces (e.g., $L^p[0, 1]$) the limit

$$\lim_{\varepsilon \downarrow 0} \varepsilon^\gamma \varphi(\varepsilon)$$

exists and is finite for the right logarithmic small ball rate $\gamma$. In this section we prove analogs of this statement for the random small ball function $\ell_\varepsilon$ in many cases.

6.1. *Sup-norm.* Recalling that $\gamma = 2$ in the case of the Wiener process and for the uniform norm on $\mathbb{C}[0, 1]$ (denoted by $\| \cdot \|$), and having in mind the deterministic equivalent of $\ell_\varepsilon$, we prove:

THEOREM 6.1. *There exists $\mathcal{K} \in (0, \infty)$ such that*

$$\lim_{\varepsilon \downarrow 0} \varepsilon^2 \ell_\varepsilon = \mathcal{K} \qquad a.s.$$

REMARK 6.2. Recall that the small ball function $\varphi$ satisfies $\lim_{\varepsilon \downarrow 0} \varepsilon^2 \varphi(\varepsilon) = \mathcal{K}_0$, where $\mathcal{K}_0 \in (0, \infty)$ is the principal eigenvalue of the Dirichlet problem on the unit ball of $\mathbb{R}^d$. Using Corollary 4.4, we can compare $\mathcal{K}$ with $\mathcal{K}_0$:

$$2\mathcal{K}_0 \leq \mathcal{K} \leq 8\mathcal{K}_0.$$

PROOF OF THEOREM 6.1. It suffices to show that the limit

(6.1)
$$\lim_{\varepsilon \downarrow 0} \varepsilon^2 \varphi_*(\varepsilon)$$

exists for $\varphi_*(\varepsilon) = \mathbb{E}[\ell_\varepsilon]$.

We slightly modify $\varphi_*$ in order to gain a transparent semi-additivity property. For $\varepsilon > 0$, let

$$(6.2) \qquad \tilde{\ell}_\varepsilon : \mathbb{C}[0, \infty) \to [0, \infty),$$

$$w \mapsto -\sup_{x \in \mathbb{R}^d} \log \mathbb{P}^x(\|W - w\| \leq \varepsilon),$$

and let $\tilde{\varphi}_*(\varepsilon) := \int \tilde{\ell}_\varepsilon \, d\mu$. Let us denote for $a \geq 0$ and $w : [0, a] \to \mathbb{R}$

$$\|w\|_{[0,a]} := \sup_{0 \leq t \leq a} |w(t)|$$

and

$$\bar{\ell}_a(w) := -\sup_{x \in \mathbb{R}^d} \log \mathbb{P}^x(\|W - w\|_{[0,a]} \leq \varepsilon).$$

Notice that $\tilde{\ell}_\varepsilon$ and $\bar{\ell}_{1/\varepsilon^2}$ are equidistributed when considering the functions as random variables on the canonical Wiener space. In particular, $\tilde{\varphi}_*(\varepsilon) = \Lambda(1/\varepsilon^2)$ for $\Lambda(a) := \int \bar{\ell}_a \, d\mu$, $a \geq 0$.

We denote by $(\theta_t)_{t \geq 0}$ the canonical ergodic flow on Wiener space, that is, for all $t \geq 0$

$$(6.3) \qquad \theta_t : \mathbb{C}[0, \infty) \to \mathbb{C}[0, \infty),$$

$$w \mapsto (\theta_t w)(s) = w(t + s) - w(s).$$

We are going to show that the family $(-\bar{\ell}_t)_{t \geq 0}$ is subadditive for the canonical dynamical system on Wiener space. Indeed, by the Markov property, one obtains for $a, b \geq 0$ that

$$
\begin{aligned}
-\bar{\ell}_{a+b}(w) \\
&= \sup_{x \in \mathbb{R}^d} \log \mathbb{P}^x(\|W - w\|_{[0,a+b]} \leq 1) \\
&= \sup_{x \in \mathbb{R}^d} \log \mathbb{P}^x(\|W - w\|_{[0,a]} \leq 1, \|W(a) - w(a) + \theta_a W - \theta_a w\|_{[0,b]} \leq 1) \\
&\leq \sup_{x \in \mathbb{R}^d} \log \mathbb{P}^x(\|W - w\|_{[0,a]} \leq 1) + \sup_{\tilde{x} \in \mathbb{R}^d} \log \mathbb{P}^{\tilde{x}}(\|W - \theta_a w\|_{[0,b]} \leq 1) \\
&= -\bar{\ell}_a(w) - \bar{\ell}_b(\theta_a w).
\end{aligned}
$$

Therefore, $\Lambda(a) = \int \bar{\ell}_a \, d\mu$ is superadditive and there exists some constant $\mathcal{K} \in [0, \infty]$ such that

$$\lim_{a \to \infty} \frac{\Lambda(a)}{a} = \mathcal{K}$$

and thus

$$(6.4) \qquad \lim_{\varepsilon \downarrow 0} \varepsilon^2 \tilde{\varphi}_*(\varepsilon) = \mathcal{K}.$$

The finiteness of $\mathcal{K}$ is easily obtained by an application of Corollary 4.4.

It remains to prove the asymptotic equivalence of $\varphi_*$ and $\tilde{\varphi}_*$. Set

$$D_\varepsilon = \{f \in \mathbb{C}[0,1] : |f(t)| \le \varepsilon, \varepsilon \le t \le 1\}$$

and consider a shift function $g_\varepsilon(t) = \min\{\varepsilon, t\}$. Then for any $x \in \mathbb{R}^d$ with $|x| \le \varepsilon$ and any $w \in \mathbb{C}[0,1]$ we have $B(w + x\mathbb{1}, \varepsilon) \subset w + \frac{g_\varepsilon}{\varepsilon}x + D_\varepsilon$. Therefore, by the Cameron–Martin formula,

$$\mu(B(w + c\mathbb{1}, \varepsilon)) \le \mu\left(w + \frac{g_\varepsilon}{\varepsilon}x + D_\varepsilon\right)$$

$$\le \mu(w + D_\varepsilon) \sup_{f \in w + D_\varepsilon} \exp\left(\frac{1}{\varepsilon}\langle g, f(\varepsilon)\rangle\right).$$

$$\le \mu(w + D_\varepsilon) \exp\left(|w(\varepsilon)| + \varepsilon\right).$$

Next, we can link $\mu(w + D_\varepsilon)$ back to conventional small ball probabilities. Indeed, it is true that

$$B(w, \varepsilon + \varepsilon^{5/4}) \supset (w + D_\varepsilon) \cap \{f : \|f_\varepsilon^0 - w_\varepsilon^0\| \le \varepsilon^{3/4}\},$$

where $w_\varepsilon^0(s) = \varepsilon^{-1/2}[w(s\varepsilon) - sw(\varepsilon)]$ and $f_\varepsilon^0$ is defined similarly by using $f$. Hence,

$$\mu(B(w, \varepsilon + \varepsilon^{5/4})) \ge \mu(w + D_\varepsilon)\mu_0(B(w_\varepsilon^0, \varepsilon^{3/4})),$$

where $\mu_0$ stands for the law of the Brownian bridge. Altogether,

$$\tilde{\ell}_\varepsilon(w) = -\sup_{c \in [-\varepsilon, \varepsilon]} \log\mu(B(w + c\mathbb{1}, \varepsilon))$$

$$\ge -\log\mu(w + D_\varepsilon) - (|w(\varepsilon)| + \varepsilon)$$

$$\ge -\log\mu(B(w, \varepsilon + \varepsilon^{5/4})) + \log\mu_0(B(w_\varepsilon^0, \varepsilon^{3/4})) - (|w(\varepsilon)| + \varepsilon).$$

From this estimate, it follows that $\varphi_* \lesssim \tilde{\varphi}_*$. On the other hand, by definition it holds $\varphi_* \ge \tilde{\varphi}_*$, and thus it follows that $\varphi_* \sim \tilde{\varphi}_*$. Now (6.4) yields the existence of the limit (6.1) and the proof is complete.  □

6.2. *Hölder norms.* We briefly discuss a modification of the previous result valid for *Hölder seminorms*. It is well known that a seminorm

$$\|f\|_\beta := \sup_{\substack{s,t \in [0,1] \\ s \ne t}} \frac{|f(t) - f(s)|}{|t - s|^\beta}$$

is $\beta$-self-similar and $\infty$-superadditive, using the terminology of [14] (see the next section). Therefore, the related small ball rate is $\gamma = (1/2 - \beta)^{-1}$ (recall that $\|W\|_\beta$ is finite iff $0 \le \beta < 1/2$). The proof of the previous section works equally well for Hölder seminorms. In the first part of the proof, the function $\Psi(a) := \tilde{\varphi}_*(a^{\beta - 1/2})$ turns out to be semi-additive. The second part of the proof is not necessary at all. Indeed, since $\|1\|_\beta = 0$, we have the identity $\varphi_* = \tilde{\varphi}_*$.

6.3. *Other norms.*  In this section, we prove the existence of small ball constants in the case of the Wiener process for a broad class of norms.

We follow the ideas of [14] and use the terminology introduced therein concerning self-similar and superadditive norms (see also [10] and [16]). Recall that a family of seminorms indexed by intervals of the real line is called *β-self-similar* iff

$$\|f(c\cdot)\|_{I/c} = c^{\beta}\|f\|_I.$$

It is called *p-superadditive* iff

$$\|f\|_{[a_0,a_n]} \geq \left(\|f\|_{[a_0,a_1]}^p + \cdots + \|f\|_{[a_{n-1},a_n]}^p\right)^{1/p} \qquad \text{for } p < +\infty,$$

$$\|f\|_{[a_0,a_n]} \geq \sup\left(\|f\|_{[a_0,a_1]}, \ldots, \|f\|_{[a_{n-1},a_n]}\right) \qquad \text{for } p = +\infty.$$

First, notice that the most interesting $\infty$-superadditive norms were considered in the two preceding sections. Therefore, in the sequel, we only consider *p*-superadditive norms with *finite p*. Again, see many examples in [14], for example, $L^p$-norms, Sobolev norms, and so on.

Let $\|\cdot\| = \|\cdot\|_{[0,1]}$ be a $\beta$-self-similar and *p*-superadditive norm. Notice that, by [14], $\gamma = (1/2 - \beta - 1/p)^{-1}$ is the right logarithmic small ball rate.

For $w \in \mathbb{C}[0,\infty)$, $\varepsilon > 0$ and $a \geq 0$ we let

$$\tilde{\ell}_{\varepsilon}(w) = -\sup_{x\in\mathbb{R}^d} \log \mathbb{P}^x(\|W-w\| \leq \varepsilon),$$

$$\Lambda_a(w) = \sup_{x\in\mathbb{R}^d} \log \mathbb{E}^x \exp\left(-\|W-w\|_{[0,a]}^p\right).$$

The functions $\tilde{\ell}_{\varepsilon}(\cdot)$ and $\Lambda_a(\cdot)$ are considered as random variables on the canonical Wiener space. We are now in a position to state the main theorem.

THEOREM 6.3.  *Assume that $\beta + 1/p < 1/2$. Then there exists a constant $\mathcal{K} \in (0,\infty)$ such that*

(6.5) $$\lim_{\varepsilon \downarrow 0} \varepsilon^{\gamma}\tilde{\ell}_{\varepsilon} = \mathcal{K} \qquad \text{in probability,}$$

*where $\gamma = (1/2 - \beta - 1/p)^{-1}$.*

REMARK 6.4.  Clearly, $\ell_{\varepsilon}$ and $\tilde{\ell}_{\varepsilon}$ are closely related. One even expects that the theorem remains true when replacing $\tilde{\ell}_{\varepsilon}$ by $\ell_{\varepsilon}$ in most cases. As we will show in the next section, we can do so if the underlying norm is the $L^p$-norm.

PROOF OF THEOREM 6.3.  The proof is based on the subadditivity of $\Lambda_a$ for the ergodic canonical flow $(\theta_t)_{t\geq 0}$. In fact, for any $w \in \mathbb{C}[0,\infty)$, $x \in \mathbb{R}^d$ and $a, b \geq 0$, by the superadditivity of our norm,

$$\log \mathbb{E}^x \exp\left(-\|W - w\|^p_{[0, a+b]}\right)$$

$$\leq \log \mathbb{E}^x \exp\left(-\|W - w\|^p_{[0,a]} - \|W - w\|^p_{[a, a+b]}\right)$$

$$= \log \mathbb{E}^x \exp\left(-\|W - w\|^p_{[0,a]} - \|W(a) - w(a) + \theta_a(W) - \theta_a(w)\|^p_{[0,b]}\right)$$

$$\leq \log \mathbb{E}^x\left[\exp\left(-\|W - w\|^p_{[0,a]}\right)\right] + \sup_{x'} \log \mathbb{E}^{x'} \exp\left(-\|W - \theta_a(w)\|^p_{[0,b]}\right).$$

Consequently, $\Lambda_{a+b}(w) \leq \Lambda_a(w) + \Lambda_b(\theta_a w)$. Subadditivity and the ergodicity imply that the following limit exists a.s.:

$$(6.6) \qquad -K := \lim_{a \to \infty} \frac{\Lambda_a}{a} \in [-\infty, 0].$$

Now, for $a \geq 0$ and $w \in \mathbb{C}[0, \infty)$, let

$$\tilde{\Lambda}_a(w) := \sup_x \log \mathbb{E}^x \exp\left(-a\|W - w\|^p_{[0,1]}\right).$$

Considered as a random variable on the canonical Wiener space, $\tilde{\Lambda}_a$ is equidistributed with $\Lambda_{a^{1/q}}$ for $q := p(1/2 - \beta)$. Hence, using (6.6),

$$\lim_{a \to \infty} a^{-1/q} \tilde{\Lambda}_a = -K \qquad \text{in probability.}$$

Now a lower bound for $\tilde{\ell}$ is obtained via the Markov inequality. For any fixed $w \in \mathbb{C}[0, 1]$, any $a, \varepsilon > 0$, $x \in \mathbb{R}^d$, we have

$$\tilde{\Lambda}_a(w) \geq \log \mathbb{E}^x \exp\{-a\|W - w\|^p\} \geq \log \mathbb{P}^x\{\|W - w\| \leq \varepsilon\} - a\varepsilon^p,$$

thus

$$\tilde{\ell}_\varepsilon(w) \geq -\tilde{\Lambda}_a(w) - a\varepsilon^p.$$

The choice of $a = a(\varepsilon) := (K/q)^{q/(q-1)} \varepsilon^{-pq/(q-1)}$ now yields

$$(6.7) \qquad \varepsilon^\gamma \tilde{\ell}_\varepsilon \gtrsim \mathcal{K} \qquad \text{in probability,}$$

where $\gamma = (1/2 - \beta - 1/p)^{-1}$ and

$$(6.8) \qquad \mathcal{K} := (q-1)(K/q)^{q/(q-1)}.$$

In particular, $\mathcal{K}$ is finite, since $\tilde{\ell}_\varepsilon$ is enclosed by $\varphi(\varepsilon)$ and $\ell_\varepsilon$ which are both of order $\varepsilon^{-\gamma}$ (see [14]).

It remains to prove the converse bound to (6.7). Toward this aim, we mimic the proof of the de Bruijn Tauberian theorem. Let $\varepsilon > 0$ and let $a = a(\varepsilon) := \frac{\mathcal{K}}{q-1} \varepsilon^{-pq/(q-1)}$. Recall that

$$(6.9) \qquad \tilde{\Lambda}_a \sim -\mathcal{K} \frac{q}{q-1} \varepsilon^{-\gamma}, \qquad \varepsilon \downarrow 0, \text{ in probability.}$$

For fixed $\varepsilon$ and $w \in \mathbb{C}[0, \infty)$ the supremum $\sup_{x \in \mathbb{R}^d} \log \mathbb{E}^x[e^{-a\|W-w\|_{[0,1]}}]$ is attained for some $x_0 = x_0(\varepsilon, w) \in \mathbb{R}^d$. Fix now $N \in \mathbb{N}$ and associate $\varepsilon$ with $\varepsilon_i := \varepsilon_i(\varepsilon) := \varepsilon i/N$, $i \in \mathbb{N}_0$. For a fixed value $\eta \in (0, 1)$, we consider

$$I := \mathbb{N} \cap [(0, (1-\eta)N + 1] \cup [(1+\eta)N, 2q^{1/p}N + 1]].$$

We estimate

$$
\begin{aligned}
\exp\{\tilde{\Lambda}_a(w)\} &= \mathbb{E}^{x_0}[\exp\{-a\|W-w\|^p\}] \\
&\leq \mathbb{E}^{x_0}[\mathbb{1}_{[(1-\eta)\varepsilon, (1+\eta)\varepsilon]}(\|W-w\|) \exp\{-a\|W-w\|^p\}] \\
&\quad + \sum_{i \in I} \mathbb{E}^{x_0}[\mathbb{1}_{[\varepsilon_{i-1}, \varepsilon_i]}(\|W-w\|) \exp\{-a\|W-w\|^p\}] \\
&\quad + \mathbb{E}^{x_0}[\mathbb{1}_{[2q^{1/p}\varepsilon, \infty)}(\|W-w\|) \exp\{-a\|W-w\|^p\}] \\
&=: \Sigma_\varepsilon^1(w) + \Sigma_\varepsilon^2(w) + \Sigma_\varepsilon^3(w).
\end{aligned}
$$

Again we consider the functions $\Sigma_\varepsilon^1$, $\Sigma_\varepsilon^2$ and $\Sigma_\varepsilon^3$ as random variables on the canonical Wiener space. We will see that $\Sigma_\varepsilon^1$ is the dominating term in the estimate. First we bound the logarithms of the summands in $\Sigma_\varepsilon^2$. One has

$$
\begin{aligned}
&\log \mathbb{E}^{x_0}[\mathbb{1}_{[\varepsilon_{i-1}, \varepsilon_i]}(\|W-w\|) \exp\{-a\|W-w\|^p\}] \\
&\quad \leq -a\varepsilon_{i-1}^p + \log \mathbb{P}^{x_0}(\|W-w\| \leq \varepsilon_i) \leq -a\varepsilon_{i-1}^p - \tilde{\ell}_{\varepsilon_i}(w).
\end{aligned}
$$

By (6.7) the right-hand side of the previous equation satisfies in probability

$$
\begin{aligned}
-a\varepsilon_{i-1}^p - \tilde{\ell}_{\varepsilon_i} &\lesssim -a\varepsilon^p(i-1)^p/N^p - \mathcal{K}\varepsilon_i^{-\gamma} \\
&= -\mathcal{K}\varepsilon^{-\gamma}\left[\frac{(i-1)^p}{(q-1)N^p} + N^\gamma/i^\gamma\right] =: -\kappa_i \mathcal{K}\varepsilon^{-\gamma}.
\end{aligned}
$$

Now let $f(x) = \frac{x^p}{(q-1)} + x^{-\gamma}$. By elementary analysis one obtains for $i \in I \cap [N+1, \infty)$

$$\kappa_i = \frac{(i-1)^p}{(q-1)N^p} + N^\gamma/i^\gamma \geq f\left(\frac{i-1}{N}\right) - \gamma N^{-1},$$

and for $i \in I \cap (0, N]$,

$$\kappa_i = \frac{(i-1)^p}{(q-1)N^p} + N^\gamma/i^\gamma \geq f\left(\frac{i}{N}\right) - \frac{p}{q-1}N^{-1}.$$

The function $f$ is strictly convex and attains its global minimum at 1. Now choose $N \in \mathbb{N}$ sufficiently large such that for all $i \in I$

$$\kappa_i > f(1) = \frac{q}{q-1}.$$

Then all summands in $\Sigma_\varepsilon^2$ are in probability of order $o(\exp\{\tilde{\Lambda}_a\})$ [see (6.9)] and one has

$$\lim_{\varepsilon\downarrow 0}\frac{\Sigma_\varepsilon^2}{\exp\{\tilde{\Lambda}_a\}}=0 \qquad \text{in probability.}$$

Moreover, $\Sigma_\varepsilon^3(w)\leq\exp\{-2^p aq\varepsilon^p\}=\exp\{-2^p\mathcal{K}\frac{q}{q-1}\varepsilon^{-\gamma}\}$ and, hence,

$$\lim_{\varepsilon\downarrow 0}\frac{\Sigma_\varepsilon^3}{\exp\{\tilde{\Lambda}_a\}}=0 \qquad \text{in probability.}$$

Therefore, $\exp\{\tilde{\Lambda}_a\}\sim\Sigma_\varepsilon^1$ in probability. Since

$$\log\Sigma_\varepsilon^1(w)\leq-a(1-\eta)^p\varepsilon^p+\log\mathbb{P}^{x_0}(\|W-w\|\leq(1+\eta)\varepsilon)$$
$$\leq-a(1-\eta)^p\varepsilon^p-\tilde{\ell}_{(1+\eta)\varepsilon}(w)$$

we arrive at

$$\tilde{\ell}_{(1+\eta)\varepsilon}\leq-a(1-\eta)^p\varepsilon^p-\log\Sigma_\varepsilon^1$$
$$\sim-a(1-\eta)^p\varepsilon^p-\tilde{\Lambda}_a$$
$$\sim\mathcal{K}\varepsilon^{-\gamma}+[1-(1-\eta)^p]\frac{\mathcal{K}}{q-1}\varepsilon^{-\gamma}.$$

Here, all equivalences hold in probability. Let $\eta\downarrow 0$, then with (6.7)

$$\lim_{\varepsilon\downarrow 0}\varepsilon^\gamma\tilde{\ell}_\varepsilon=\mathcal{K} \qquad \text{in probability.} \qquad\qquad \square$$

REMARK 6.5.  It would be very interesting to extend the results of this subsection to self-similar processes other than the Wiener process, as done in [14] for fractional Brownian motion in the case of nonrandomly centered balls. Our main subadditivity argument for upper bound seems to fail in the non-Markovian case.

6.4.  $L^p$-norm.  In the sequel $\|\cdot\|$ denotes $L^p[0,1]$-norm for some fixed $p\in[1,\infty)$. Moreover, let $\|f\|_{[a,b]}$ and $\|f\|_{[a,b],\infty}$ $(a\leq b)$ denote the $L^p$-norm and the sup-norm over the interval $[a,b]$, respectively. For $\varepsilon>0$ we consider the map

$$\ell_\varepsilon:\mathbb{C}[0,\infty)\to[0,\infty),$$
$$w\mapsto-\log\mathbb{P}(\|W-w\|\leq\varepsilon)$$

as random variable on the canonical Wiener space.

Our objective is to prove:

THEOREM 6.6. *One has*

$$\lim_{\varepsilon \downarrow 0} \varepsilon^2 \ell_\varepsilon = \mathcal{K} \qquad a.s., \tag{6.10}$$

*where $\mathcal{K}$ is as in Theorem 6.3.*

Notice that, in order to prove the theorem, it suffices to prove convergence (6.10) in probability. Since clearly $\tilde{\ell}_\varepsilon \lesssim \ell_\varepsilon$, it remains to show that $\ell_\varepsilon \lesssim \mathcal{K}\varepsilon^{-2}$ in probability.

We need some preliminary propositions.

PROPOSITION 6.7. *For $w \in \mathbb{C}[0,1]$ and $\varepsilon \in (0,1/2)$ it is true that*

$$\sup_{z \in \mathbb{R}^d} \mathbb{P}^z(\|W - w\| \le \varepsilon) \le \sup_{t \in [\varepsilon, 2\varepsilon]} \sup_{x \in B(w_t, \varepsilon^{1-1/p})} \mathcal{F}(x, t, \varepsilon),$$

*where $\mathcal{F}(x, t, \varepsilon') = \mathbb{P}^{x,t}(\|W - w\|_{[t,1]} \le \varepsilon')$ for $x \in \mathbb{R}^d, t \in [0,1]$ and $\varepsilon' \in \mathbb{R}$.*

PROOF. Note that if $\|W - w\| \le \varepsilon$, then the stopping time

$$T := \inf\{t \ge \varepsilon : |W_t - w_t| \le \varepsilon^{1-1/p}\}$$

satisfies $T \le 2\varepsilon$. Using the Markov property of the Wiener process, we obtain

$$\begin{aligned}
\mathbb{P}^z(\|W - w\| \le \varepsilon) &= \mathbb{E}^z[\mathbb{1}_{\{T \le 2\varepsilon\}} \mathcal{F}(W_T, T, (\varepsilon^p - \|W - w\|_{[0,T]}^p)^{1/p})] \\
&\le \mathbb{E}^z[\mathbb{1}_{\{T \le 2\varepsilon\}} \mathcal{F}(W_T, T, \varepsilon)] \\
&\le \sup_{t \in [\varepsilon, 2\varepsilon]} \sup_{x \in B(w_t, \varepsilon^{1-1/p})} \mathcal{F}(x, t, \varepsilon). \qquad \square
\end{aligned}$$

PROPOSITION 6.8. *Let $\varepsilon, \theta \in (0, 1/2)$ and $w \in \mathbb{C}[0,1]$. Then*

$$\mathbb{P}(\|W - w\|_{[0,2\varepsilon],\infty} \le \varepsilon) \cdot \sup_{t \in [\varepsilon, 2\varepsilon]} \sup_{x \in B(w_t, \varepsilon^{1-1/p})} \mathcal{F}(x, t, \varepsilon)$$

$$\le (\theta/3)^{-d} \mathbb{P}(\|W - w\|_{[0,1]} \le \varepsilon(1 + \tilde{\theta})) \exp\left(2\theta^{-p}\varepsilon^{-1}[3 + 2\|w\|_{[0,2\varepsilon],\infty}]\right),$$

*where $\tilde{\theta} = (2\varepsilon)^{1/p} + 5\theta$.*

PROOF. We fix $t \in [\varepsilon, 2\varepsilon]$ and $x \in B(w_t, \varepsilon^{1-1/p})$. Let us consider

$$A = \{f \in \mathbb{C}[0, \infty) : \|f - w\|_{[0,t],\infty} \le \varepsilon\}, \qquad A_y = \{f \in A : f(t) \in B(y, \theta\varepsilon)\}$$

and choose $y \in B(w(t), \varepsilon)$ such that

$$\mu(A_y) \ge (\theta/3)^d \mu(A). \tag{6.11}$$

We have

$$|x - y| \le |x - w(t)| + |y - w(t)| \le \varepsilon^{1-1/p} + \varepsilon \le 2. \tag{6.12}$$

Next, define a shift function $g \in H$ by $g'(s) = t^{-1}\theta^{-p}\mathbb{1}_{[(1-\theta^p)t,t]}$. Obviously, $g(t) = 1$ and $\|g\|_{[0,t]} \le \theta t^{1/p}$. Let $Q_{x,y} = A_y + (x-y)g$.

For any $h \in Q_{x,y}$ we have two properties: $h(t) \in B(x,\theta\varepsilon)$ and, using (6.12),

$$\|h - w\|_{[0,t]} \le \varepsilon t^{1/p} + |x - y|\|g\|_{[0,t]}$$
$$\le \varepsilon t^{1/p} + 2\varepsilon^{1-1/p}\theta t^{1/p}$$
$$\le \varepsilon(t^{1/p} + 4\theta).$$

We also need an elementary inequality

$$\inf_{z \in B(x,\theta\varepsilon)} \mathbb{P}^{t,z}(\|W - w\|_{[t,1]} \le \varepsilon(1+\theta)) \ge \mathcal{F}(x,t,\varepsilon).$$

It follows that

$$\mathbb{P}(\|W - w\|_{[0,1]} \le \varepsilon(1 + t^{1/p} + 5\theta))$$
$$\ge \mathbb{P}(\|W - w\|_{[0,t]} \le \varepsilon(t^{1/p} + 4\theta); \|W - w\|_{[t,1]} \le \varepsilon(1+\theta))$$
$$\ge \mathbb{P}(W \in Q_{x,y}) \inf_{z \in B(x,\theta\varepsilon)} \mathbb{P}^{t,z}(\|W - w\|_{[t,1]} \le \varepsilon(1+\theta))$$
$$\ge \mathbb{P}(W \in Q_{x,y})\mathcal{F}(x,t,\varepsilon).$$

Now we pass from $\mathbb{P}(W \in Q_{x,y})$ to $\mathbb{P}(W \in A_y)$. Recall that $Q_{x,y} = A_y + (x-y)g$. Hence, by the Cameron–Martin formula

$$\mathbb{P}(W \in Q_{x,y}) \ge \mathbb{P}(W \in A_y)\exp\{-\mathcal{M}\},$$

where

$$\mathcal{M} := \tfrac{1}{2}|x-y|^2|g|_\mu^2 + |x-y|t^{-1}\theta^{-p}\sup_{h \in A_y}\left|\int_{(1-\theta^p)t}^t dh(s)\right|.$$

By using (6.12) and the definition of $g$ we have

$$\tfrac{1}{2}|x-y|^2|g|_\mu^2 \le 2t^{-1}\theta^{-p}.$$

Moreover, for $h \in A_y \subset A$,

$$\left|\int_{(1-\theta^p)t}^t dh(s)\right| = |h(t) - h((1-\theta^p)t)| \le 2(\varepsilon + \|w\|_{[0,t],\infty}).$$

Hence,

$$\mathcal{M} \le 2\theta^{-p}\varepsilon^{-1}(1 + 2\varepsilon + 2\|w\|_{[0,t],\infty})$$

and by combining this bound with (6.11) and other previous estimates we get

$$\mathbb{P}(\|W - w\|_{[0,1]} \le \varepsilon(1 + t^{1/p} + 5\theta))$$
$$\ge (\theta/3)^d\mathbb{P}(W \in A)\exp\left\{-2\theta^{-p}\varepsilon^{-1}(3 + 2\|w\|_{[0,t],\infty})\right\}\mathcal{F}(x,t,\varepsilon),$$

and thus the result follows. □

PROOF OF THEOREM 6.6. Recall that by Theorem 6.3:

$$\tilde{\ell}_\varepsilon \sim \mathcal{K} \frac{1}{\varepsilon^2} \qquad \text{as } \varepsilon \downarrow 0, \text{ in probability.}$$

Moreover, by the above proposition one has for $w \in \mathbb{C}[0,\infty)$, $\varepsilon, \theta \in (0, 1/2)$,

$$\ell_{\varepsilon(1+\tilde{\theta})}(w) \leq \tilde{\ell}_\varepsilon(w) - \log \mathbb{P}(\|W - w\|_{[0,2\varepsilon],\infty} \leq \varepsilon)$$
$$- d \log(\theta/3) + 2\theta^{-p}\varepsilon^{-1}[3 + 2\|w\|_{[0,2\varepsilon],\infty}],$$

where $\tilde{\theta}$ is as in Proposition 6.8. Choosing $\theta := \varepsilon^{1/(2p)}$ we obtain

$$\ell_{\varepsilon(1+7\varepsilon^{1/(2p)})}(w) \leq \tilde{\ell}_\varepsilon(w) - \log \mathbb{P}(\|W - w\|_{[0,2\varepsilon],\infty} \leq \varepsilon)$$
$$- d \log(\varepsilon^{1/(2p)}/3) + 2\varepsilon^{-3/2}[3 + 2\|w\|_{[0,2\varepsilon],\infty}].$$

Now let $w$ be a $\mu$-distributed random variable. Then all summands but $\ell_{\varepsilon(1+7\varepsilon^{1/(2p)})}(w)$ and $\tilde{\ell}_\varepsilon(w)$ are of order $o(\varepsilon^{-2})$ in probability. Consequently,

$$\ell_{\varepsilon(1+7\varepsilon^{1/(2p)})} \lesssim \tilde{\ell}_\varepsilon \sim \mathcal{K} \frac{1}{\varepsilon^2} \qquad \text{in probability.}$$

The assertion follows when choosing $\tilde{\varepsilon} > 0$ with $\tilde{\varepsilon} = \varepsilon(1 + 7\varepsilon^{1/(2p)})$ and letting $\tilde{\varepsilon} \downarrow 0$. □

Fakultät II
Institut für Mathematik
MA 7-5
Technische Universität Berlin
Strasse des 17. Juni 136
10623 Berlin
Germany
e-mail: dereich@math.tu-berlin.de

Faculty of Mathematics and Mechanics
St. Petersburg State University
Bibliotechnaya pl. 2
198504 Stary Peterhof
Russia
e-mail: lifts@mail.rcom.ru